# HOLOMORPHIC RATIONAL FUNCTIONS OF SEVERAL VARIABLES AND SUMS OF SQUARES OF POLYNOMIALS

M.F. Bessmertnyi

Necessary and sufficient conditions are obtained under which the numerator of the partial derivative of a rational function holomorphic in an open upper poly-halfplane is the sum of squares of polynomials.

## 1. Introduction

If a polynomial $p$ in $d$ variables with real coefficient can be written as a sum of squares of polynomials, then $p$ must take only nonnegative values in $\mathbb{R}^d$. Not every non-negative polynomial of $d$ variables is the sum of squares of polynomials. Hilbert's 17th problem on the representation of a non-negative polynomial as a sum of squares of rational functions was solved by Artin in 1927. There are quite a few articles on this topic (see, for example, the bibliography in reviews [11], [13])

If a rational function $p/q$ is holomorphic in an open upper poly-halfplane, then non-negative polynomials can appear as partial Wronskians

$$W_k[q,p] = q(z)\frac{\partial p(z)}{\partial z_k} - p(z)\frac{\partial q(z)}{\partial z_k} \tag{1.1}$$

of the pair of polynomials such that $q(z) \neq 0$, $p(z) \neq 0$ for $z$ from an open upper poly-halfplane [1 - 3, 9, 10, 17]. If a non-negative polynomial cannot be represented as a sum of squares of polynomials (examples of Motzkin, Robinson [15], M.-D. Choi [6]), then the analysis shows it also cannot be represented in the form (1.1).

The problem arises: *for a rational function $f = p/q$ holomorphic in an open upper poly-halfplane, find the necessary and sufficient conditions for the representation of the partial Wronskian $W_{z_k}[q,p]$ as a sum of squares of polynomials.*

In this article, we prove the following statement.

**Main Theorem**. *Let $f(z) = p(z)/q(z)$ be a rational function with real coefficients holomorphic on open upper poly-halfplane. The partial Wronskian $W_{z_1}[q,p]$ is a sum of squares of polynomials if and only if for each fixed $x_2,\cdots,x_d \in \mathbb{R}$*

$$\operatorname{Im}\frac{p(z_1,x_2,\cdots,x_d)}{q(z_1,x_2,\cdots,x_d)} \geq 0, \qquad \operatorname{Im} z_1 > 0.$$

The paper is organized as follows. In Section 2, we explain the terminology and provide preliminary information.

Section 3 studies the properties of denominators of rational functions in the Artin decomposition into the sum of squares. Theorem 3.3 and Proposition 3.4 allows localizing the singularities of rational functions in Artin's decomposition.

A convenient representation for the partial Wronskians $W_{z_k}[q,p]$ is given in Theorem 4.1 (Product Polarization Theorem). In fact, for a rational function $f = p/q$, this theorem implies Hefer's expansion $f(z)-f(\varsigma)=\sum(z_k-\varsigma_k)F_k(z,\varsigma)$ with additional conditions of symmetry $F_k(z,\varsigma)=\Psi(z)B_k\Psi(\varsigma)^T$, where $B_k$ are symmetric real matrices.

In Section 5, we study the set of Gram matrices of a given polynomial of degree $2n$ and prove the Representation Defect Lemma. This lemma describes the uncertainty of the matrix pencil from Theorem 4.1.

A representation of a rational function with one nonnegative partial Wronskian is obtained in Theorem 6.1. The representation contains the Artin denominator of the nonnegative partial Wronskian in explicit form.

In Section 7, we prove the Main Theorem (Theorem 7.4).

## 2. Preliminaries

$\mathcal{H}^d=\{z\in\mathbb{C}^d : \operatorname{Im} z_1>0,\cdots,\operatorname{Im} z_d>0\}$ is an open upper poly-halfplane.

$\mathbb{R}[z]=\mathbb{R}[z_1,\cdots,z_d]$ is a ring of polynomials with real coefficients.

$\mathcal{Z}_V(I)$ is an zero set (in $V=\mathbb{R}$ or $\mathbb{C}$) of the family $I$ of functions.

$\mathcal{I}(S)$ is an ideal of polynomials in $\mathbb{R}[z]$ vanishing on $S$. The principal ideal generated by the polynomial $p$ is denoted by $(p)$.

An ideal $I$ in $\mathbb{R}[z]$ is said to *be real* if, for every sequence $a_1,\cdots,a_p$ of elements of $\mathbb{R}[z]$, we have

$$a_1^2+...+a_p^2\in I \quad\Rightarrow\quad a_i\in I,\ \text{ for } i=1,\cdots,p.$$

**Theorem 2.1** ([4], Theorem 4.5.1). *Let $\mathcal{R}$ be real closed field and $p$ an irreducible polynomial in $\mathcal{R}[z_1,\cdots,z_d]$. Then the following properties are equivalent:*

(i) *The ideal $(p)$ is real.*

(ii) $(p)=\mathcal{I}(\mathcal{Z}_{\mathcal{R}}(p))$.

(iii) *The polynomial $p$ has a nonsingular zero in $\mathcal{R}^d$ (i.e., there is an $x\in\mathcal{R}^d$ such that $p(x)=0$ and $\dfrac{\partial p}{\partial z_i}(x)\neq 0$ for some $i\in(1,\cdots,d)$.*

(iv) *The sign of the polynomial $p$ changes on $\mathcal{R}^d$ (i.e., $p(x)p(y)<0$ for some $x,y$ in $\mathcal{R}^d$).*

(v) $\dim(\mathcal{Z}_{\mathcal{R}}(p))=d-1$.

If $\mathcal{R}$ is a field, then $\mathcal{R}(x_1,\cdots,x_d)$ denotes the set of rational functions in variables $x_1,\cdots,x_d$ with coefficients from the field $\mathcal{R}$.

**Artin's Theorem.** ([12], Ch. XI, Corollary 3.3). *Let $\mathcal{R}$ be a real field admitting only one ordering. Let $f(x)\in\mathcal{R}(x_1,\cdots,x_d)$ be a rational function that does not take negative values: $f(a)\geq 0$ for all $a=(a_1,\cdots,a_d)\in\mathcal{R}^d$, in which $f(a)$ is defined. Then $f(x)$ there is the sum of squares in $\mathcal{R}(x_1,\cdots,x_d)$.*

Any matrix $A$ is called *real*, if $\bar{A} = A$ (where the bar denotes the replacement of each element of $A$ by a complex conjugate number). The symbol $A^T$ denotes the transposition operation. If $A$ is a matrix with complex elements, then $A^* = \bar{A}^T$ is the Hermitian conjugate matrix.

Any real symmetric $m \times m$ matrix $A$ is called *a positive semidefinite* $(A \geq 0)$ if the inequality $\eta^T A \eta \geq 0$ holds for all $\eta \in \mathbb{R}^m$, and *a positive definite* $(A > 0)$ if $\eta^T A \eta > 0$ for all $\eta \neq 0$.

A polynomial $F(z)$ will be called *a PSD polynomial* (positive semidefinite) if $F(x) \geq 0$ for all $x \in \mathbb{R}^d$. A PSD polynomial $F(z)$ will be called *an SOS polynomial* (sum of squares) if $F(z) = H(z)H(z)^T$, where $H(z)$ is a row vector of polynomials.

If $F(z)$ is an SOS polynomial, then $s(z)^2 F(z)$ is also an SOS for each polynomial $s(z)$. If $F(z)$ is not representable as a sum of squares of polynomials, then the question arises: *for which* $s(z)$ *is the polynomial* $s(z)^2 F(z)$ *also not an SOS polynomial?*

Any polynomial $s(z)$ in $\mathbb{R}[z_1, \cdots, z_d]$ is said to *be indefinite* if, for some $x, y$ in $\mathbb{R}^d$ we have $s(x)s(y) < 0$.

**Proposition 2.2.** ([7], Lemma 2.1). *Let* $F(x)$ *be a PSD not SOS polynomial and* $s(x)$ *an irreducible indefinite polynomial in* $\mathbb{R}[x_1, \cdots, x_d]$. *Then* $s^2 F$ *is also a PSD not SOS polynomial.*

*Proof.* Clearly $s^2 F$ is PSD. If $s^2 F = \sum_k h_k^2$, then for every real tuple $a$ with $s(a) = 0$, it follows that $s^2 F(a) = 0$. This implies $h_k(a)^2 = 0 \;\; \forall k$. So on the real variety $s = 0$, we have $h_k = 0$ as well. By Theorem 2.1 for each $k$, there exists $g_k$ so that $h_k = s g_k$.This gives $F = \sum_k g_k^2$, which is a contradiction. □

If $\alpha = (\delta_1, \cdots, \delta_d) \in \mathbb{N}_0^d$ is an integer multi-index, then $z^\alpha = z_1^{\delta_1} \cdots z_d^{\delta_d}$ is a monomial. Let $\{z^{\alpha_j}\}_{j=1}^M$ be a set of all monomials of degree $\leq n$ in variables $z_1, \cdots, z_d$. Each polynomial $F(z)$ of degree $2n$ can be represented as

$$F(z) = (z^{\alpha_1}, \cdots, z^{\alpha_M}) \begin{pmatrix} a_{11} & \cdots & a_{1M} \\ \vdots & \ddots & \vdots \\ a_{1M} & \cdots & a_{MM} \end{pmatrix} \begin{pmatrix} z^{\alpha_1} \\ \vdots \\ z^{\alpha_M} \end{pmatrix}. \tag{2.1}$$

The symmetric matrix $A = \{a_{jk}\}_{j,k=1}^M$ is called *a Gram matrix* of $F(z)$. The Gram matrix is not uniquely determined by $F(z)$. It is known ([14, Theorem 1]) that PSD polynomial $F(z)$ is an SOS if and only if $F(z)$ has a positive semidefinite Gram matrix.

## 3. Zero Sets of Artin's Denominator

Let $F(z) \in \mathbb{R}[z_1, \cdots, z_d]$ be a PSD not SOS polynomial. By Artin's theorem, there exists a polynomial $s(z)$ such that $s(z)^2 F(z)$ is an SOS polynomial. The polynomial $s(z)$ is called *Artin's denominator of* $F(z)$.

**Proposition 3.1.** *Suppose* $s(z)^2 F(z)$ *is an SOS polynomial and each irreducible factor of* $s(z)$ *is indefinite; then* $F(z)$ *is also an SOS polynomial.*

*Proof*. Suppose $F(z)$ is a PSD not SOS polynomial. Let $s(z) = s_1(z)\cdots s_m(z)$ be the decomposition of $s(z)$ into irreducible factors. Successively applying Proposition 2.2 to the polynomials

$$F_1(z) = s_1^2(z)F(z),\ F_2(z) = s_2^2(z)F_1(z),\ \cdots\ ,\ F_m(z) = s_m^2(z)F_{m-1}(z),$$

we obtain $F_m(z) = s^2(z)F(z)$ is a PSD not SOS polynomial. Contradiction. □

**Definition 3.2.** Artin's denominator $s(z)$ of a PSD not SOS polynomial $F(z)$ is called *a minimal Artin denominator* if polynomial $s(z)/s_j(z)$ is not Artin's denominator of $F(z)$ for each of the irreducible factor $s_j(z)$ of $s(z)$.

**Theorem 3.3.** *Each PSD not SOS polynomial $F(z)$ has a non-constant minimum Artin denominator $s(z)$. All irreducible factors of $s(z)$ do not change sign on $\mathbb{R}^d$.*

*Proof*. By Artin's theorem, there exists a polynomial $r(z)$ for which $r(z)^2F(z)$ is an SOS. Each irreducible factor of $r(z)$ is either indefinite or does not change sign on $\mathbb{R}^d$. Then $r(z) = r_0(z)s(z)$, where all irreducible factors of $s(z)$ do not change sign on $\mathbb{R}^d$, and irreducible factors of $r_0(z)$ are indefinite. Consider the polynomial $F_1(z) = s(z)^2F(z)$. By condition, $r_0(z)^2F_1(z) = r(z)^2F(z)$ is an SOS. The irreducible factors of $r_0(z)$ are indefinite. By Proposition 3.1, $F_1 = s^2F$ is an SOS. Then $s(z)$ is also Artin's denominator for $F(z)$. Let $s_0$ be some irreducible factor of $s$. If $s/s_0$ remains the Artin denominator of $F$, then the factor $s_0$ is removed from $s$. Removing all "excess" irreducible factors from the polynomial $s(z)$, we obtain Artin's denominator with the required properties. □

**Proposition 3.4**. *Let $s(z) \in \mathbb{R}[z_1,\cdots,z_d]$ be a irreducible polynomial that does not change sign on $\mathbb{R}^d$. If $\partial s(z)/\partial z_1 \not\equiv 0$, then there exists a point $z' = (z_1', x_2', \cdots, x_d')$, $\operatorname{Im} z_1' > 0$, $x_2',\cdots,x_d' \in \mathbb{R}$, for which $s(z') = 0$. Moreover, there exists $z'' \in \mathcal{H}^d$ for which $s(z'') = 0$.*

*Proof*. For $d = 1$ this is true. Let us $d \geq 2$. Suppose $s(z_1, x_2, \cdots, x_d)$ has only real zeros in $z_1$ for all $x_2,\cdots,x_d \in \mathbb{R}$. Then the equation $s(x) = 0$ defines a real manifold of dimension $d-1$. By Theorem 2.1, $(s)$ is a real ideal and irreducible polynomial $s(z)$ is indefinite. This contradicts the condition. Complex zeros of $s_x(z_1) = s(z_1, x_2, \cdots, x_d)$ form complex conjugate pairs. Then there exists $x_2',\cdots,x_d' \in \mathbb{R}$ and $z_1' \in \mathbb{C}$, $\operatorname{Im} z_1' > 0$ such that $s(z_1', x_2', \cdots, x_d') = 0$. Put $z_j'' = x_j' + iy_j$, $j = 2,\cdots,d$. If $y_j > 0$ sufficiently small, then the polynomial $s(z_1, z_2'', \cdots, z_d'')$ still vanishes at some point $z_1''$ from the open upper half-plane. □

## 4. Product Polarization Theorem

Recall if $\alpha = (\delta_1,\cdots,\delta_d) \in \mathbb{N}_0^d$ is a multi-index, then $z^\alpha = z_1^{\delta_1}\cdots z_d^{\delta_d}$ is a monomial.

The following theorem gives a special representation for the product of two polynomials.

**Theorem 4.1. (Product Polarization Theorem)**. *Let $q(z)$, $p(z)$, $z \in \mathbb{C}^d$ be real polynomials for which $\max\{\deg q,\ \deg p\} = n$, $\max\{\deg_{z_k} q,\ \deg_{z_k} p\} = n_k$, $k = 1,\cdots,d$. If* $\Psi(z)=(z^{\alpha_1},\cdots,z^{\alpha_N})$ *is a row vector of all monomial satisfying the conditions* $\deg z^{\alpha_j}\le n$, $\deg_{z_k} z^{\alpha_j}\le n_k$, $k=1,\cdots,d$ *, then there exist real symmetric matrices* $A_k$, $k=0,1,\cdots,d$ *such that*

$$q(\varsigma)p(z)=\Psi(\varsigma)(A_0+z_1A_1+...+z_dA_d)\,\Psi(z)^T, \qquad \varsigma,z\in\mathbb{C}^d, \tag{4.1}$$

$$q(z)\partial p(z)/\partial z_k-p(z)\partial q(z)/\partial z_k=\Psi(z)A_k\Psi(z)^T,\ k=1,\cdots,d. \tag{4.2}$$

We need some lemmas.

**Lemma 4.2**. *Let* $k\ge 0$ *be an integer and* $\varsigma^{\mu_1}=\varsigma_2\varsigma_4\cdots\varsigma_{2k}$, $\varsigma^{\nu}=\varsigma_1\varsigma_3\cdots\varsigma_{2k+1}$. *Then there exist real symmetric* $(2k+1)\times(2k+1)$ *matrices* $C_k$, $k=1,\cdots,2k+1$ *and multiaffine monomials* $\{\varsigma^{\mu_j}\}_{j=2}^{2k+1}$ *of degree* $k$ *in variables* $\varsigma_1,\cdots,\varsigma_{2k+1}$ *such that*

$$(\varsigma_1C_1+...+\varsigma_{2k+1}C_{2k+1})\begin{pmatrix}\varsigma^{\mu_1}\\ \varsigma^{\mu_2}\\ \vdots\\ \varsigma^{\mu_{2k+1}}\end{pmatrix}=\begin{pmatrix}\varsigma^{\nu}\\ 0\\ \vdots\\ 0\end{pmatrix}. \tag{4.3}$$

*Proof*. If $k=0$, then $\varsigma^{\mu_1}=1$ (empty product) and $\varsigma^{\nu}=\varsigma_1$. Then $\varsigma^{\nu}=\varsigma_1\varsigma^{\mu_1}$, and the matrix pencil $C(\varsigma)=\varsigma_1\cdot 1$ has a size of $1\times 1$. If $k\ge 1$, then the multiaffine monomials $\{\varsigma^{\mu_j}\}_{j=2}^{2k+1}$ are defined by the relations

$$\varsigma^{\mu_2}=\varsigma^{\nu}/\varsigma_1=\varsigma_3\varsigma_5\cdots\varsigma_{2k+1},\quad \varsigma^{\mu_j}=\varsigma_{j-2}\varsigma^{\mu_{j-2}}/\varsigma_{j-1},\qquad j=3,4,\cdots,2k+1.$$

Note that $\varsigma^{\mu_{2k+1}}=\varsigma^{\nu}/\varsigma_{2k+1}$, $\varsigma^{\mu_{2k}}=\varsigma_{2k+1}\varsigma^{\mu_1}/\varsigma_{2k}$. Let us define the matrix pencil $C(\varsigma)=\{c_{ij}(\varsigma)\}_{i,j=1}^{2k+1}$:

$$c_{ij}(\varsigma)=\begin{cases}(-1)^{\max\{i,j\}}\varsigma_{\min\{i,j\}}/2, & |i-j|=1\\ \varsigma_{\max\{i,j\}}/2, & |i-j|=2k\\ 0, & \text{otherwise}\end{cases}.$$

It is easy to see that $\overline{C(\bar z)}=C(z)=C(z)^T$. Let us calculate the components $b_i(\varsigma)$, $i=1,\cdots,2k+1$ of the right-hand side of (4.3).

$$b_1(\varsigma)=\sum_{j=1}^{2k+1}c_{1j}(\varsigma)\varsigma^{\mu_j}=c_{12}(\varsigma)\varsigma^{\mu_2}+c_{1,2k+1}(\varsigma)\varsigma^{\mu_{2k+1}}=\frac{1}{2}\varsigma_1\frac{\varsigma^{\nu}}{\varsigma_1}+\frac{1}{2}\varsigma_{2k+1}\frac{\varsigma^{\nu}}{\varsigma_{2k+1}}=\varsigma^{\nu}.$$

For $2\le i\le 2k$ we get

$$b_i(\varsigma)=\sum_{j=1}^{2k+1}c_{ij}(\varsigma)\varsigma^{\mu_j}=c_{i,i-1}(\varsigma)\varsigma^{\mu_{i-1}}+c_{i,i+1}(\varsigma)\varsigma^{\mu_{i+1}}=\frac{(-1)^i}{2}\varsigma_{i-1}\varsigma^{\mu_{i-1}}+\frac{(-1)^{i+1}}{2}\varsigma_i\varsigma^{\mu_{i+1}}$$

$$=\left[\varsigma^{\mu_{i+1}}=\varsigma_{i-1}\varsigma^{\mu_{i-1}}/\varsigma_i\right]=\left((-1)^i\varsigma_{i-1}\varsigma^{\mu_{i-1}}+(-1)^{i+1}\varsigma_{i-1}\varsigma^{\mu_{i-1}}\right)/2=0.$$

For $i=2k+1$ we obtain

$$b_{2k+1}(\varsigma)=\sum_{j=1}^{2k+1}c_{2k+1,j}(\varsigma)\varsigma^{\mu_j}=c_{2k+1,1}(\varsigma)\varsigma^{\mu_1}+c_{2k+1,2k}(\varsigma)\varsigma^{\mu_{2k}}=\frac{1}{2}\left(\varsigma_{2k+1}\varsigma^{\mu_1}-\varsigma_{2k}\varsigma^{\mu_1}\frac{\varsigma_{2k+1}}{\varsigma_{2k}}\right)=0.\ \square$$

**Lemma 4.3**. *Let* $z^{\alpha_1}$, $z^{\beta}$, $z=(z_1,\cdots,z_d)\in\mathbb{C}^d$ *be monomials satisfy the conditions*

$$\max\{\deg z^{\alpha_1},\ \deg z^{\beta}\}=n,\quad \max\{\deg_{z_k} z^{\alpha_1},\ \deg_{z_k} z^{\beta}\}=n_k,\ k=1,\cdots,d.$$

*Then there exist matrices* $B_k=\bar{B}_k=B_k^T$, $k=0,1,\cdots,d$ *such that*

$$(B_0+z_1B_1+...+z_dB_d)\begin{pmatrix} z^{\alpha_1}\\ z^{\alpha_2}\\ \vdots\\ z^{\alpha_N}\end{pmatrix}=\begin{pmatrix} z^{\beta}\\ 0\\ \vdots\\ 0\end{pmatrix}, \tag{4.4}$$

*where* $\{z^{\alpha_j}\}_{j=1}^N$ *are all monomials satisfying* $\deg z^{\alpha_j}\le n$, $\deg_{z_k} z^{\alpha_j}\le n_k$, $k=1,\cdots,d$.

*Proof*. We multiply the monomials $z^{\alpha_1}$ and $z^{\beta}$ by some powers of the new variable $z_0$ so that

$$z^{\alpha_1}\mapsto z_0^{l_\alpha}z^{\alpha_1},\ z^{\beta}\mapsto z_0^{l_\beta}z^{\beta}\quad\text{and}\quad \deg z_0^{l_\alpha}z^{\alpha_1}=m,\quad \deg z_0^{l_\beta}z^{\beta}=m+1.$$

Let $\hat{z}^{\gamma}$, $\hat{z}=(z_0,z_1,\cdots,z_d)$ be the greatest common divisor of the monomials $z_0^{l_\alpha}z^{\alpha_1}$, $z_0^{l_\beta}z^{\beta}$. Then $z_0^{l_\alpha}z^{\alpha_1}=m_1(\hat{z})\hat{z}^{\gamma}$, $z_0^{l_\beta}z^{\beta}=m_2(\hat{z})\hat{z}^{\gamma}$, where the subsets of the variables of the monomials $m_1(\hat{z})$ and $m_2(\hat{z})$ do not intersect. If $\deg m_1(\hat{z})=k\ge 0$, then $\deg m_2(\hat{z})=k+1$. By Lemma 4.2, there exist matrices $C_j=\bar{C}_j=C_j^T$, $j=1,\cdots,2k+1$ such that

$$(\varsigma_1C_1+...+\varsigma_{2k+1}C_{2k+1})\begin{pmatrix} \varsigma_2\varsigma_4\cdots\varsigma_{2k}\hat{z}^{\gamma}\\ \varsigma^{\mu_2}\hat{z}^{\gamma}\\ \vdots\\ \varsigma^{\mu_{2k+1}}\hat{z}^{\gamma}\end{pmatrix}=\begin{pmatrix} \varsigma_1\varsigma_3\cdots\varsigma_{2k+1}\hat{z}^{\gamma}\\ 0\\ \vdots\\ 0\end{pmatrix}, \tag{4.5}$$

where $\{\varsigma^{\mu_j}\}_{j=2}^{2k+1}$ are multiaffine monomials of degree $k$ in variables $\varsigma_1,\varsigma_2,\cdots,\varsigma_{2k+1}$.

In (4.5) we replace the variables $\varsigma_2,\varsigma_4,\cdots,\varsigma_{2k}$ by the variables of the monomial $m_1(\hat{z})$, and the variables $\varsigma_1,\varsigma_3,\cdots,\varsigma_{2k+1}$ by the variables of the monomial $m_2(\hat{z})$, so that

$$\varsigma_2\varsigma_4\cdots\varsigma_{2k}\hat{z}^{\gamma}\ \mapsto\ m_1(\hat{z})\hat{z}^{\gamma}=z_0^{l_\alpha}z^{\alpha_1},\quad \varsigma_1\varsigma_3\cdots\varsigma_{2k+1}\hat{z}^{\gamma}\ \mapsto\ m_2(\hat{z})\hat{z}^{\gamma}=z_0^{l_\beta}z^{\beta}.$$

From (4.5), setting $z_0=1$, we obtain

$$(D_{j_0}+z_{j_1}D_{j_1}+...+z_{j_r}D_{j_r})\begin{pmatrix} z^{\hat{\alpha}_1}\\ z^{\hat{\alpha}_2}\\ \vdots\\ z^{\hat{\alpha}_{2k+1}}\end{pmatrix}=\begin{pmatrix} z^{\beta}\\ 0\\ \vdots\\ 0\end{pmatrix},\qquad z^{\hat{\alpha}_1}=z^{\alpha_1}, \tag{4.6}$$

where $z_{j_1},\cdots,z_{j_r}$ are the variables of monomials $m_1(\hat{z})$, $m_2(\hat{z})$ different from the variable $z_0$. The matrices $D_{j_s}$ are the sums of the corresponding matrices $C_j$ from (4.5).

Since each monomial $\varsigma^{\mu_j}$, $j=2,\cdots,2k+1$, contains only a part of the variables $\varsigma_i$ with even indices, we see that for the variables $z_i \in m_1(\hat{z})\hat{z}^{\gamma}$ holds $\deg_{z_i} z^{\hat{\alpha}_j} \le \deg_{z_i} z^{\alpha_1} \le n_i$. Similarly, for $z_v \in m_2(\hat{z})\hat{z}^{\gamma}$ we get $\deg_{z_v} z^{\hat{\alpha}_j} \le \deg_{z_v} z^{\beta} \le n_v$.

Let $z^{\alpha_1}, z^{\alpha_2}, \cdots, z^{\alpha_M}$ be pairwise distinct monomials from the set $z^{\hat{\alpha}_1}, z^{\hat{\alpha}_2}, \cdots, z^{\hat{\alpha}_{2k+1}}$. Since $z^{\alpha_1} = z^{\hat{\alpha}_1}$, we see that there exists a $(2k+1)\times M$ matrix $B$ such that

$$\begin{pmatrix} z^{\hat{\alpha}_1} \\ z^{\hat{\alpha}_2} \\ \vdots \\ z^{\hat{\alpha}_{2k+1}} \end{pmatrix} = \begin{pmatrix} 1 & 0 & \cdots & 0 \\ b_{21} & b_{22} & \cdots & b_{2M} \\ \vdots & \vdots & \vdots & \vdots \\ b_{2k+1,1} & b_{2k+1,2} & \cdots & b_{2k+1,M} \end{pmatrix} \begin{pmatrix} z^{\alpha_1} \\ z^{\alpha_2} \\ \vdots \\ z^{\alpha_M} \end{pmatrix}.$$

From (4.6) we obtain

$$(B_0 + z_{j_1} B_{j_1} + \ldots + z_{j_r} B_{j_r}) \begin{pmatrix} z^{\alpha_1} \\ z^{\alpha_2} \\ \vdots \\ z^{\alpha_M} \end{pmatrix} = \begin{pmatrix} z^{\beta} \\ 0_2 \\ \vdots \\ 0_M \end{pmatrix}, \tag{4.7}$$

where $(B_0 + z_{j_1} B_{j_1} + \ldots + z_{j_r} B_{j_r}) = B^T (D_{j_0} + z_{j_1} D_{j_1} + \ldots + z_{j_r} D_{j_r}) B$. We extend the set $\{z^{\alpha_j}\}_{j=1}^{M}$ to the set $\{z^{\alpha_i}\}_{i=1}^{N}$ of all pairwise distinct monomials satisfying conditions $\deg z^{\alpha_i} \le n$, $\deg_{z_k} z^{\alpha_i} \le n_k$, $k=1,\cdots,d$. Supplementing the matrices in (4.7) with zero entries, we obtain (4.4). □

*Proof of Theorem 4.1.* Suppose $q(z) = \sum_{j=1}^{N} a_j z^{\alpha_j}$, $p(z) = \sum_{v=1}^{l} b_v z^{\beta_v}$. By Lemma 4.3, for fixed monomials $z^{\alpha_j}$, $z^{\beta_v}$ there exists a symmetric real matrix pencil $B_{jv}(z)$ for which

$$B_{jv}(z) \begin{pmatrix} z^{\alpha_1} \\ \vdots \\ z^{\alpha_j} \\ \vdots \\ z^{\alpha_N} \end{pmatrix} = \begin{pmatrix} 0 \\ \vdots \\ z^{\beta_v} \\ \vdots \\ 0 \end{pmatrix} - j\text{-}th\ row, \quad v=1,\cdots,l, \quad j=1,\cdots,N, \tag{4.8}$$

where $\deg z^{\alpha_i} \le n$, $\deg_{z_k} z^{\alpha_i} \le n_k$, $k=1,\cdots,d$.

We define $A(z) = A_0 + z_1 A_1 + \ldots + z_d A_d = \sum_{j=1}^{N} a_j \sum_{v=1}^{l} b_k B_{jv}(z)$. Then

$$A(z) \begin{pmatrix} z^{\alpha_1} \\ \vdots \\ z^{\alpha_N} \end{pmatrix} = \begin{pmatrix} a_1 p(z) \\ \vdots \\ a_N p(z) \end{pmatrix}. \tag{4.9}$$

Since $a_j, b_v \in \mathbb{R}$ and $\overline{B_{jv}(\overline{z})} = B_{jv}(z) = B_{jv}(z)^T$, we see that $\overline{A(\overline{z})} = A(z) = A(z)^T$. Multiplying (4.9) on the left by the row vector $(\varsigma^{\alpha_1}, \cdots, \varsigma^{\alpha_N})$, $\varsigma \in \mathbb{C}^d$, we obtain (4.1). Relations (4.2) follow from (4.1). □

## 5. Uncertainty of Representing Matrix Pencil

Let $h(z)$ be a polynomial such that $\deg h(z) = n$, $\deg_{z_k} h(z) \le n_k$, $k = 1,\cdots,d$. Then

$$H(z_0, z) = z_0^n h(z_1/z_0,\cdots,z_d/z_0)$$

is the $n$ - form (homogeneous polynomial of degree $n$) with the following properties:

$$\deg H(z_0,z) = n, \quad \deg_{z_0} H(z_0,z) \le n, \quad \deg_{z_k} H(z_0,z) \le n_k, \; k = 1,\cdots,d.$$

The polynomial $h(z)$ is PSD (SOS) if and only if the form $H(z_0,z)$ is PSD (respectively SOS). Therefore, when studying Gram`s matrices of a PSD polynomials, it suffices to restrict ourselves to the study of the Gram matrices of $2n$ - forms.

Let $\{z^{\alpha_j}\}_{i,j=1}^{N}$ be a set of all monomials of degree $n$ satisfying the conditions $\deg_{z_k} z^{\alpha_v} \le n_k$, $k = 1,\cdots,d$ and let $F(z)$ be a real $2n$-form such that $\deg_{z_k} F(z) \le 2n_k$, $k = 1,\cdots,d$. Suppose $F(z)$ has two Gram matrices $A_1$, $A_2$: $F(z) = \Psi(z)A_1\Psi(z)^T = \Psi(z)A_2\Psi(z)^T$, where $\Psi(z) = (z^{\alpha_1},\cdots,z^{\alpha_N})$. The symmetric real matrix $S = A_1 - A_2 = \{s_{ij}\}_{i,j=1}^{N}$ satisfies the relation

$$\Psi(z)S\Psi(z)^T = \sum_{i,j=1}^{N} s_{ij} z^{\alpha_i} z^{\alpha_j} = 0. \tag{5.1}$$

The set of matrices $S$ satisfying the condition (5.1) is a linear space $L_0$. Now we construct a special basis in this linear space.

**Proposition 5.1.** *Let $L_0$ be a linear space of real symmetric matrices satisfying condition* (5.1). *Then there exists a basis in $L_0$ such that the nonzero submatrices of the basis matrices are located at the intersection of rows and columns corresponding to monomials of a special form:*

$$(z_r^2 z^\gamma,\; z_r z_l z^\gamma,\; z_l^2 z^\gamma)\begin{pmatrix} 0 & 0 & -1 \\ 0 & 2 & 0 \\ -1 & 0 & 0 \end{pmatrix}\begin{pmatrix} z_r^2 z^\gamma \\ z_r z_l z^\gamma \\ z_l^2 z^\gamma \end{pmatrix} \equiv 0, \tag{5.2}$$

$$(z_r z^{\gamma_1},\; z_l z^{\gamma_1},\; z_l z^{\gamma_2},\; z_r z^{\gamma_2})\begin{pmatrix} 0 & 0 & 1 & 0 \\ 0 & 0 & 0 & -1 \\ 1 & 0 & 0 & 0 \\ 0 & -1 & 0 & 0 \end{pmatrix}\begin{pmatrix} z_r z^{\gamma_1} \\ z_l z^{\gamma_1} \\ z_l z^{\gamma_2} \\ z_r z^{\gamma_2} \end{pmatrix} \equiv 0, \quad z^{\gamma_1} \ne z^{\gamma_2}. \tag{5.3}$$

**Remark 5.2.** *Since* $\sum_{i,j=1}^{N} s_{ij} z^{\alpha_i} z^{\alpha_j} = \sum_k c_k z^{\beta_k} = 0$, and $z^{\beta_i} \ne z^{\beta_j}$, $i \ne j$, *we see that*

$$c_k = \sum_{\alpha_i + \alpha_j = \beta_k} s_{ij} = 0. \tag{5.4}$$

*If the sum* (5.4) *contains $m \ge 2$ different elements $s_{ij}$, then $m-1$ elements can be chosen as arbitrary. Then multi-index $\beta_k$ defines a $(m-1)$-dimensional subspace in $L_0$.*

We need some lemmas.

Let $\beta=(r_1,\cdots,r_d)$ be a multi-index and let $\Pi_\beta$ be a set of all unordered pairs $\pi=(z^{\alpha_i},z^{\alpha_j})$ such that $z^{\alpha_i}z^{\alpha_j}=z^\beta$, where $\deg_{z_k} z^{\alpha_\nu}\le n_k$ for all $\alpha_\nu$. It is easy to see that $z_1^{\delta_1}\cdots z_d^{\delta_d}\in\pi$ iff

$$\max\{(r_k-n_k),0\}\le\delta_k\le\min\{r_k,n_k\},\quad k=1,\cdots,d. \tag{5.5}$$

**Definition 5.3.** If $\delta_r<n_r$, $\delta_l>0$, then a map

$$z_1^{\delta_1}\cdots z_r^{\delta_r}\cdots z_l^{\delta_l}\cdots z_d^{\delta_d}\mapsto z_1^{\delta_1}\cdots z_r^{\delta_r+1}\cdots z_l^{\delta_l-1}\cdots z_d^{\delta_d} \tag{5.6}$$

is called *an elementary monomial transformation*. Any element $\pi_2\in\Pi_\beta$ is called *an elementary transformation of an element* $\pi_1\in\Pi_\beta$, if some monomial in $\pi_2$ is an elementary transformation of a monomial in $\pi_1$.

**Lemma 5.4.** *For any $\pi_1,\pi_2\in\Pi_\beta$, $\pi_1\ne\pi_2$ there exists a connecting $\pi_1$ and $\pi_2$ chain of elements $\pi_k\in\Pi_\beta$ such that each next element is an elementary transformation of the previous one.*

*Proof*. Suppose $\beta=(r_1,\cdots,r_d)$. We have $z^{\alpha_i}=z_1^{\delta_1}\cdots z_d^{\delta_d}\in\pi_1$, $z^{\alpha_j}=z_1^{\nu_1}\cdots z_d^{\nu_d}\in\pi_2$, $z^{\alpha_i}\ne z^{\alpha_j}$. If $k_s=\nu_s-\delta_s$, $s=1,\cdots,d$, then $\sum_{s=1}^d k_s=0$, $-\min\{r_s,n_s\}\le k_s\le\min\{r_s,n_s\}$. The tuple $(k_1,\cdots,k_d)$ is a componentwise sum of the elementary tuples. Each elementary tuple $(m_1,\cdots,m_d)$ contains only two nonzero components $+1$, $-1$. The tuple $(m_1,\cdots,m_d)$ corresponds elementary monomial transformation (5.6). Starting from $z^{\alpha_i}\in\pi_1$, at each step we get a monomial satisfying (5.5). At the last step, we get $z^{\alpha_j}\in\pi_2$. The element $\pi_k=(z^{\alpha_{k_1}},z^{\alpha_{k_2}})\in\Pi_\beta$ is uniquely determined by one of the monomials $z^{\alpha_{k_1}},z^{\alpha_{k_2}}$. □

**Lemma 5.5.** *If $\Pi_\beta$ contain $m\ge2$ elements, then in the set $\Pi_\beta$ there exist $(m-1)$ different pairs $\{\pi_i,\pi_j\}$ such that $\pi_j$ is an elementary transformation of $\pi_i$.*

*Proof*. Let us associate the finite graph with the set $\Pi_\beta$. The vertices are elements of the set $\Pi_\beta$. The edges form pairs $\{\pi_i,\pi_j\}$ of elements connected by an elementary transformation. By Lemma 5.4, the graph is connected. Then the graph tree contains $m-1$ edges. In the graph tree different edges are incident to different pairs of vertices. □

*Proof of Proposition 5.1*. The linear space $L_0$ is the direct sum of subspaces $L_{\beta_k}$, each of which corresponds to its own multi-index $\beta_k$. Let us construct a basis in each of these subspaces.

Suppose $\Pi_{\beta_k}$ contain $m\ge2$ elements. By Lemma 5.5, in $\Pi_{\beta_k}$ there exist $m-1$ different pairs $\{\pi_i,\pi_j\}$ such that $\pi_i$ and $\pi_j$ are connected by an elementary transformation. Let us show that each such pair $\{\pi_i,\pi_j\}$ defines a basis matrix of the form (5.2) or (5.3). The following cases are possible:

(a). One of the elements of a pair $\{\pi_1,\pi_2\}$ has the form $\pi_1=(z^{\alpha_i},z^{\alpha_i})$, and $\pi_2=(z^{\alpha_j},z^{\alpha_l})$. Then there exist a variables $z_r,z_l$ such that $z^{\alpha_j}=z_r^2z^\gamma$, $z^{\alpha_i}=z_rz_lz^\gamma$, $z^{\alpha_l}=z_l^2z^\gamma$. This triplet of monomials corresponds to the basis matrix (5.2).

(b). $\pi_1=(z^{\alpha_i},z^{\alpha_j})$, $z^{\alpha_i}\neq z^{\alpha_j}$; $\pi_2=(z^{\alpha_l},z^{\alpha_s})$, $z^{\alpha_l}\neq z^{\alpha_s}$, and monomials $z^{\alpha_i}$, $z^{\alpha_l}$ are connected by an elementary transformation. Then there exist variables $z_r,z_l$ such that $z^{\alpha_i}=z_r z^{\gamma_1}$, $z^{\alpha_l}=z_l z^{\gamma_1}$, $z^{\alpha_j}=z_l z^{\gamma_2}$, $z^{\alpha_s}=z_r z^{\gamma_2}$. Note that $z^{\gamma_1}\neq z^{\gamma_2}$. Indeed, if $z^{\gamma_1}=z^{\gamma_2}$, then $\pi_1=\pi_2$, which is impossible. The constructed four monomials define the basis matrix (5.3).

All pairs $\{\pi_i,\pi_j\}$ are different. Then corresponding $(m-1)$ matrices are linearly independent.□

**Lemma 5.6. (Representation Defect Lemma)**. *Let* $\Psi(z)=(z^{\alpha_1},\cdots,z^{\alpha_N})$ *be a row vector of all monomials satisfy the conditions* $\deg z^{\alpha_i}\le n$, $\deg_{z_k} z^{\alpha_i}\le n_k$, $k=1,\cdots,d$. *If for a real symmetric* $N\times N$ *- matrix* $S_d$ *the identities hold*

$$\Psi(z)S_d\Psi(z)^T\equiv 0,\qquad S_d\left(\partial^{n_d}\Psi(z)^T/\partial z^{n_d}\right)\equiv 0,$$

*then there exist the real symmetric* $N\times N$ *- matrices* $S_k$, $k=0,1,\cdots,d-1$ *for which*

$$(S_0+z_1S_1+...+z_{d-1}S_{d-1}+z_dS_d)\Psi(z)^T\equiv 0. \tag{5.7}$$

*Proof*. It suffices to consider the case of homogeneous monomials of degree $n$ in variables $z_0,z_1,\cdots,z_d$ and at the end put $z_0=1$.

Without loss of generality, we can assume that $\Psi(z)=(z_d^{n_d}\varphi(\hat z),\ \psi(\hat z,z_d))$, where $\hat z=(z_0,z_1,\cdots,z_{d-1})$, $\deg_{z_d}\psi(\hat z,z_d)\le n_d-1$ and $\deg_{z_d} z_d^{n_d}\varphi(\hat z)=n_d$. Then

$$S_d\ \frac{\partial^{n_d}\Psi(z)^T}{\partial z^{n_d}}=\begin{pmatrix} S_{11} & S_{12}\\ S_{12}^T & \hat S_d\end{pmatrix}\begin{pmatrix} n_d!\varphi(\hat z)^T\\ 0\end{pmatrix}\equiv 0.$$

Hence $S_{11}=0$, $S_{12}^T=0$ and $\psi(\hat z,z_d)\hat S_d\ \psi(\hat z,z_d)^T\equiv 0$. We rewrite (5.7) in block form

$$\begin{pmatrix} S_{11}(\hat z) & S_{12}(\hat z)\\ S_{12}(\hat z)^T & z_d\hat S_d+S_{22}(\hat z)\end{pmatrix}\begin{pmatrix} z_d^{n_d}\varphi(z)^T\\ \psi(\hat z,z_d)^T\end{pmatrix}=\begin{pmatrix}0\\0\end{pmatrix}. \tag{5.8}$$

Let us find a solution $\{S_{ij}(\hat z)\}_{i,j=1}^2$ of the equation (5.8) for the matrices $\hat S_{d,j}$ of the form (5.2), (5.3). For the basis matrix (5.2) there exist monomials $z_dz_rz^\gamma$ and $z_dz_lz^\gamma$ such that

$$\begin{pmatrix} 0 & 0 & 0 & -z_l & z_r\\ 0 & 0 & z_l & -z_r & 0\\ 0 & z_l & 0 & 0 & -z_d\\ -z_l & -z_r & 0 & 2z_d & 0\\ z_r & 0 & -z_d & 0 & 0\end{pmatrix}\begin{pmatrix} z_dz_rz^\gamma\\ z_dz_lz^\gamma\\ z_r^2z^\gamma\\ z_rz_lz^\gamma\\ z_l^2z^\gamma\end{pmatrix}\equiv 0.$$

Similarly, for basis matrix (5.3) and monomials $z_dz^{\gamma_1}$, $z_dz^{\gamma_2}$ we have

$$\begin{pmatrix} 0 & 0 & -z_v & z_k & 0 & 0 \\ 0 & 0 & 0 & 0 & -z_k & z_v \\ -z_v & 0 & 0 & 0 & z_d & 0 \\ z_k & 0 & 0 & 0 & 0 & -z_d \\ 0 & -z_k & z_d & 0 & 0 & 0 \\ 0 & z_v & 0 & -z_d & 0 & 0 \end{pmatrix} \begin{pmatrix} z_d z^{\gamma_2} \\ z_d z^{\gamma_1} \\ z_k z^{\gamma_1} \\ z_v z^{\gamma_1} \\ z_v z^{\gamma_2} \\ z_k z^{\gamma_2} \end{pmatrix} \equiv 0 .$$

Matrix $\hat{S}_d$ is a linear combination of basis matrices $\hat{S}_{1,j}$ of the form (5.2), (5.3).Then solution of equation (5.8) is a linear combination of such solutions. Assuming $z_0 = 1$ in the resulting solution, we obtain (5.7). □

## 6. Rational Functions with PSD not SOS Wronskian

Let $f(z) = p(z)/q(z)$ be a rational function. If the partial Wronskian $W_{z_k}[q, p]$ is a PSD polynomial, then the following statement holds.

**Theorem 6.1.** *Let* $f(z) = p(z)/q(z)$, $z \in \mathbb{C}^d$ *be rational function with real coefficients and* $\max\{\deg p, \deg q\} = n$, $\max\{\deg_{z_1} p, \deg_{z_1} q\} = n_1$. *Let* $s(z)^2 W_{z_1}[q, p] = H(z)H(z)^T$ *be an SOS polynomial for some polynomial* $s(z)$. *Then there exist a real symmetric matrices* $A_0, A_1, \cdots, A_d$, *where* $A_1$ *is positive semidefinite, such that*

$$f(z) = \frac{p(z)}{q(z)} = \frac{\Psi(\varsigma)}{q(\varsigma)s(\varsigma)}(A_0 + z_1 A_1 + ... + z_d A_d)\frac{\Psi(z)^T}{q(z)s(z)}, \qquad \varsigma, z \in \mathbb{C}^d , \tag{6.1}$$

$$W_{z_k}[q, p] = \frac{\Psi(z)}{s(z)} A_k \frac{\Psi(z)^T}{s(z)}, \qquad k = 1, \cdots, d , \tag{6.2}$$

$$\Psi(z) A_1 \Psi(z)^T = H(z)H(z)^T , \tag{6.3}$$

*where* $\Psi(z) = (z^{\alpha_1}, \cdots, z^{\alpha_N})$ *is a row vector of all monomial* $z^{\alpha_i}$ *satisfies the conditions* $\deg z^{\alpha_i} \le n + \deg s(z) = m$, $\deg_{z_1} z^{\alpha_i} \le n_1 + \deg_{z_1} s(z) = m_1$.

*Proof.* By Theorem 4.1, there exists a matrix pencil $B_0 + z_1 B_1 + ... + z_d B_d$ for which

$$q(\varsigma)s(\varsigma)p(z)s(z) = \Psi(\varsigma)(B_0 + z_1 B_1 + ... + z_d B_d)\Psi(z)^T , \tag{6.4}$$

$$W_{z_k}[qs, ps] = \Psi(z) B_k \Psi(z)^T , \quad k = 1, \cdots, d , \tag{6.5}$$

where a row vector $\Psi(z) = (z^{\alpha_1}, \cdots, z^{\alpha_N})$ satisfies the conditions of the theorem. Differentiating (6.4) $(m_1 + 1)$ times in $z_1$, we obtain

$$B_1 \frac{\partial^{m_1} \Psi(z)^T}{\partial z_1^{m_1}} \equiv 0 . \tag{6.6}$$

By condition, $s(z)^2 W_{z_1}[q, p]$ is an SOS form. Then there exists a matrix $A_1 \ge 0$ for which

$$\Psi(z) B_1 \Psi(z)^T = W_{z_1}[qs, ps] = s^2 W_{z_1}[q, p] = H(z)H(z)^T = \Psi(z) A_1 \Psi(z)^T .$$

Since $\deg_{z_1} W_{z_1}[qs, ps] \le 2m_1 - 2$, we see that if $\deg_{z_1} a_{ii} z^{\alpha_i} z^{\alpha_i} = 2m_1$, then $a_{ii} = 0$. From $A_1 \ge 0$ we get

$$A_1 \frac{\partial^{m_1} \Psi(z)^T}{\partial z_1^{m_1}} \equiv 0 . \tag{6.7}$$

By (6.7) and (6.6) it follows $S_1 = A_1 - B_1$, $S_1^T = S_1 = \bar{S}_1$ satisfies the conditions Representation Defect Lemma. Then there exist real symmetric matrices $S_0, S_2, \cdots, S_d$ such that

$$\Psi(\varsigma)(S_0 + z_1 S_1 + z_2 S_2 + ... + z_d S_d)\Psi(z)^T \equiv 0, \qquad \varsigma, z \in \mathbb{C}^d . \tag{6.8}$$

Adding (6.8) to (6.4) and dividing both sides of the resulting identity by $q(\varsigma)s(\varsigma)q(z)s(z)$, we obtain (6.1). Relations (6.2), (6.3) follow from identities $s^2 W_{z_k}[q, p] = W_{z_k}[qs, ps] = \Psi(z) A_k \Psi(z)^T$. □

## 7. Main Theorem

In the univariate case, coprime polynomials have no common zeros. For several variables, the situation is different (the simplest example $z_1 / z_2$).

Let $Z_{\mathbb{C}}(h) = \{z \in \mathbb{C}^d : h(z) = 0\}$ be a zero set of the polynomial $h(z) \in \mathbb{R}[z_1, \cdots, z_d]$.

We need an analogue of Theorem 1.3.2 from [16]:

**Proposition 7.1**. Let $s(z)$, $h(z)$ *be coprime polynomials in the ring* $\mathbb{R}[z_1, \cdots, z_d]$ *and* $s(z_1', \hat{z}') = h(z_1', \hat{z}') = 0$ *for fixed* $z_1' \in \mathbb{C}$, $\hat{z}' \in \mathcal{R}^{d-1}$, *where* $\mathcal{R} = \mathbb{R}$ *or* $\mathbb{C}$. *Put* $\Omega = \Omega_1 \times \Omega_{\mathcal{R}}$, *where* $\Omega_1 \subseteq \mathbb{C}$ *is a some neighborhood of* $z_1'$, *and* $\Omega_{\mathcal{R}} \subseteq \mathcal{R}^{d-1}$ *is a some neighborhood of* $\hat{z}'$. *Then neither of the sets* $Z(s) \cap \Omega$ *and* $Z(h) \cap \Omega$ *is a subset of the other.*

*Proof*. Without loss of generality, we can assume $s(z_1, \hat{z}') \not\equiv 0$ and $h(z_1, \hat{z}') \not\equiv 0$. If this is not the case, then we perform a linear transformation of variables $z_k \mapsto z_1 + z_k$, $k = 2, \cdots, d$. There exists a neighborhood $\Delta' \subseteq \mathcal{R}^{d-1}$ of $\hat{z}'$ such that each $\hat{z} \in \Delta'$ corresponds to at least one point $(z_1, \hat{z}) \in \mathcal{Z}_{\mathbb{C}}(s) \cap \Omega$.

The factorial ring $\mathcal{K} = \mathbb{R}[z_2, \cdots, z_d]$ contains no zero divisors. Let $\mathcal{F}$ be the field of fractions of ring $\mathcal{K}$. We will consider $s(z_1, \hat{z})$ and $h(z_1, \hat{z})$ as elements of the ring $\mathcal{K}[z_1]$.The polynomials $s(z_1, \hat{z})$ and $h(z_1, \hat{z})$ are coprime in the ring $\mathcal{K}[z_1]$. By Gauss's lemma (see [8]), $s(z_1, \hat{z})$ and $h(z_1, \hat{z})$ are coprime in the ring $\mathcal{F}[z_1]$. Then there exist polynomials $A, B \in \mathcal{F}[z_1]$ for which

$$As + Bh = 1. \tag{7.1}$$

Let $C(\hat{z}) \in \mathcal{K} = \mathbb{R}[z_2, \cdots, z_d]$ be the common denominator of the all coefficients of the polynomials $A$ and $B$. Multiplying (7.1) by $C(\hat{z})$, we obtain

$$\tilde{A}(z_1, \hat{z}) s(z_1, \hat{z}) + \tilde{B}(z_1, \hat{z}) h(z_1, \hat{z}) = C(\hat{z}), \tag{7.2}$$

where $\tilde{A}(z_1, \hat{z})$ and $\tilde{B}(z_1, \hat{z})$ are polynomials. Since $C(\hat{z})$ is a nonzero polynomial, we see that $C(\hat{z}') \ne 0$ for some $\hat{z}' \in \Delta' \subseteq \mathcal{R}^{d-1}$. Then there is $z_1'$ such that $(z_1', \hat{z}') \in \mathcal{Z}_{\mathbb{C}}(s) \cap \Omega$. From (7.2) it follows $h(z_1, \hat{z}) \ne 0$. □

**Lemma 7.2**. *Let* $f(z)=p(z)/q(z)$ *be a rational function holomorphic on* $\mathcal{H}^d=\{z\in\mathbb{C}^d:\operatorname{Im}z_1>0,\cdots,\operatorname{Im}z_d>0\}$. *Then a function* $f_x(\varsigma)=f(\varsigma,x_2,\cdots,x_d)$ *is holomorphic on* $\mathcal{H}^1$ *for any fixed* $x_2',\cdots,x_d'\in\mathbb{R}$.

*Proof*. If $q(\varsigma',x_2',\cdots,x_d')=0$, $\operatorname{Im}\varsigma'>0$, $x_2',\cdots,x_d'\in\mathbb{R}$, then for a sufficiently small positive $y_2'>0,\cdots,y_d'>0$ the polynomial $q(\varsigma,x_2'+iy_2',\cdots,x_d'+iy_d')$ still vanishes at the point from the open upper half-plane. Contradiction. □

**Lemma 7.3**. *Let* $p(z)$, $q(z)$ *be polynomials with real coefficients. If partial Wronskian* $W_{z_1}[q,p]$ *is a PSD polynomial, then for each fixed* $x_2,\cdots,x_d\in\mathbb{R}$

$$\operatorname{Im}\frac{p(z_1,x_2,\cdots,x_d)}{q(z_1,x_2,\cdots,x_d)}\ge 0,\qquad \operatorname{Im}z_1>0.$$

*Proof*. By Artin`s Theorem, there exists polynomial $s(z)$ for which

$$s(z)^2W_{z_1}[q,p]=s(z)^2\left(q(z)\frac{\partial p(z)}{\partial z_1}-p(z)\frac{\partial q(z)}{\partial z_1}\right)$$

is an SOS polynomial. By Theorem 6.1, there exist a real symmetric matrices $A_0,A_1,\cdots,A_d$, where $A_1$ is positive semidefinite, such that

$$f(z)=\frac{p(z)}{q(z)}=\frac{\Psi(\varsigma)}{q(\varsigma)s(\varsigma)}(A_0+z_1A_1+...+z_dA_d)\frac{\Psi(z)^T}{q(z)s(z)},\ \varsigma,z\in\mathbb{C}^d,$$

where $\Psi(z)$ is row vector of monomials. This implies

$$\operatorname{Im}f(z_1,\hat{x})=\operatorname{Im}z_1\frac{\Psi(z_1,\hat{x})}{q(z_1,\hat{x})s(z_1,\hat{x})}A_1\frac{\Psi(z_1,\hat{x})^*}{\overline{q(z_1,\hat{x})s(z_1,\hat{x})}}\ge 0,\ \operatorname{Im}z_1>0.$$ □

**Theorem 7.4. (Main Theorem)**. *Let* $f(z)=p(z)/q(z)\in\mathbb{R}(z_1,\cdots,z_d)$ *be a function holomorphic on* $\mathcal{H}^d$. *The partial Wronskian* $W_{z_1}[q,p]$ *is an SOS polynomial if and only if for each fixed* $x_2,\cdots,x_d\in\mathbb{R}$

$$\operatorname{Im}\frac{p(z_1,x_2,\cdots,x_d)}{q(z_1,x_2,\cdots,x_d)}\ge 0,\qquad \operatorname{Im}z_1>0. \tag{7.3}$$

*Proof*. Since each SOS polynomial is PSD, we see that (7.3) follows from Lemma 7.1.

Let us the converse. Let $W_{z_1}[q,p]$ be a PSD not SOS polynomial and $s(z)$ be its minimal Artin's denominator. Then $s^2(z)W_{z_1}[q,p]=H(z)H(z)^T$, where $H(z)$ is row vector of polynomials. Each irreducible factor $s_0(z)$ of $s(z)$ cannot be a divisor of all elements from $H(z)$, otherwise, $\hat{s}=s/s_0$ is also Artin's denominator of $W_{z_1}[q,p]$, which contradicts the minimality of $s(z)$.

By Theorem 6.1, there exists a symmetric matrix pencil $A(z)=A_0+z_1A_1+...+z_dA_d$ with a positive semidefinite matrix $A_1\ge 0$ such that

$$f(z)=\frac{p(z)}{q(z)}=\frac{\Psi(\varsigma)}{q(\varsigma)s(\varsigma)}(A_0+z_1A_1+...+z_dA_d)\frac{\Psi(z)^T}{q(z)s(z)}, \tag{7.4}$$

$$W_{z_1}[q,p]=\frac{\Psi(z)}{s(z)}A_1\,\frac{\Psi(z)^T}{s(z)}=\frac{H(z)}{s(z)}\,\frac{H(z)^T}{s(z)}\,, \tag{7.5}$$

where $\Psi(z)$ is a row vector of monomials.

From (7.4) we get

$$\operatorname{Im} f(z) \;=\; \operatorname{Im} z_1 \frac{H(z)}{q(z)s(z)}\,\overline{\frac{H(z)^*}{q(z)s(z)}} \;+\; \sum_{k=2}^{d}\operatorname{Im} z_k\,\frac{\Psi(z)}{q(z)s(z)}A_k\,\overline{\frac{\Psi(z)^*}{q(z)s(z)}}\,. \tag{7.6}$$

There are two possibilities for the irreducible factor $s_0(z)$: (**a**) $\partial s_0(z)/\partial z_1 \not\equiv 0$, (**b**) $\partial s_0(z)/\partial z_1 \equiv 0$.

***Case*** **(a)**. $\partial s_0(z)/\partial z_1 \not\equiv 0$. By propositions 3.4 and 7.1, there exists a point $z'=(z_1',\hat{x}')$, $\operatorname{Im} z_1'>0$, $\hat{x}'\in\mathbb{R}^{d-1}$ such that $s(z')=s_0(z')=0$, $h(z')\neq 0$. From (7.6) we get

$$\lim_{\varsigma\to z_1'}\operatorname{Im} f(\varsigma,\hat{x}') \;=\; \lim_{\varsigma\to z_1'}\operatorname{Im}\varsigma\cdot\sum_j\frac{\left|h_j(\varsigma,\hat{x})\right|^2}{\left|q(\varsigma,\hat{x})\right|^2\left|s(\varsigma,\hat{x})\right|^2}=+\infty\,, \tag{7.7}$$

which contradicts the condition (7.3).

***Case*** **(b)**. If $\partial s_0(z)/\partial z_1 \equiv 0$, then for some $k\neq 1$ $\partial s_0(z)/\partial z_k \not\equiv 0$ holds. By proposition 3.4, there exists a point $z''\in\mathcal{H}^{d-1}$ such that $s_0(z'')=0$. Then for each $z_1\in\mathbb{C}$, except for a finite number of points, the inequality $h(z_1,z'')\neq 0$ holds.

Let $\Omega\subseteq\mathcal{H}^{d-1}$ be a neighborhood of the point $z''\in\mathcal{H}^{d-1}$. The function $f(x_1',z)$ is holomorphic on $\mathcal{H}^{d-1}$, if fixed $x_1'\in\mathbb{R}$. Then $\operatorname{Im} f(x_1',z)$ is bounded for $z\in\Omega$. From (7.6) we get

$$\operatorname{Im} f(x_1',z)=\sum_{k=2}^{d}\operatorname{Im} z_k\,\frac{\Psi(x_1',z)}{q(x_1',z)s(z)}A_k\,\overline{\frac{\Psi(x_1',z)^*}{q(x_1',z)s(z)}}\,. \tag{7.8}$$

The sum

$$\sum_{k=2}^{d}\operatorname{Im} z_k\,\frac{\Psi(z_1,z)}{q(z_1,z)s(z)}A_k\,\overline{\frac{\Psi(z_1,z)^*}{q(z_1,z)s(z)}}$$

is a rational function of variables $z_1$, $z\in\mathbb{C}^{d-1}$, taking uniformly bounded values at $z_1=x_1'\in\mathbb{R}$, $z\in\Omega$. Then there exists a neighborhood $\Delta'\subseteq\mathbb{C}$ of the point $x_1'$ such that for all $z_1\in\Delta'$ and $z\in\Omega$ holds

$$\left|\sum_{k=2}^{d}\operatorname{Im} z_k\,\frac{\Psi(z_1,z)}{q(z_1,z)s(z)}A_k\,\overline{\frac{\Psi(z_1,z)^*}{q(z_1,z)s(z)}}\right|\le M<+\infty\,. \tag{7.9}$$

Since for each $z_1\in\mathbb{C}$, except for a finite number of points, the inequality $h(z_1,z'')\neq 0$ holds, we see that there exists $z_1''\in\Delta'$, $\operatorname{Im} z_1''>0$ such that $h(z_1'',z'')\neq 0$. From (7.6) we get

$$\operatorname{Im} f(z_1,z'') \;=\; \operatorname{Im} z_1\cdot\sum_j\frac{\left|h_j(z_1,z'')\right|^2}{\left|q(z_1,z'')\right|^2\left|s(z'')\right|^2} \;+\; \sum_{k=2}^{d}\operatorname{Im} z_k\,\frac{\Psi(z_1,z'')}{q(z_1,z'')s(z'')}A_k\,\overline{\frac{\Psi(z_1,z'')^*}{q(z_1,z'')s(z'')}}\,. \tag{7.10}$$

According to (7.9), the second term in (7.10) is bounded for $z_1 \in \Delta'$. The first term increases indefinitely at $z_1 \to z_1''$. Hence $\lim_{z_1 \to z_1''} \operatorname{Im} f_{ii}(z_1, z'') = +\infty$. Then the function $f(z_1, z)$ has a singularity at the point $(z_1'', z'') \in \mathcal{H}^d$. Contradiction. □